\documentclass[12pt]{amsart}
\usepackage[margin=1in]{geometry}
\usepackage{amsthm}
\usepackage{amsmath,amssymb,xypic}
\usepackage[colorlinks,citecolor=blue]{hyperref}
\usepackage{epigraph}
\usepackage{graphicx}
\usepackage{natbib}
\usepackage{fancyhdr} 
\usepackage{enumitem}
\usepackage{mathrsfs}
\newcommand{\be}{\begin{equation}}
\newcommand{\ee}{\end{equation}}
\newcommand{\bes}{\begin{equation*}}
\newcommand{\ees}{\end{equation*}}
\newcommand{\bea}{\begin{eqnarray}}
\newcommand{\eea}{\end{eqnarray}}
\newcommand{\beas}{\begin{eqnarray}}
\newcommand{\eeas}{\end{eqnarray}}
\newcommand{\ben}{\begin{note}}
\newcommand{\een}{\end{note}}
\newcommand{\bexl}{\vskip0.1em\noindent\hrulefill\vskip1em\begin{ExerciseList}}
\newcommand{\eexl}{\end{ExerciseList}\hrulefill}

\newcommand{\bthm}{\begin{theorem}}
\newcommand{\ethm}{\end{theorem}}
\newcommand{\bpro}{\begin{prop}}
\newcommand{\epro}{\end{prop}}
\newcommand{\bcor}{\begin{corollary}}
\newcommand{\ecor}{\end{corollary}}
\newcommand{\bcon}{\begin{conjecture}}
\newcommand{\econ}{\end{conjecture}}
\newcommand{\bp}{\begin{proof}}
\newcommand{\ep}{\end{proof}}
\newcommand{\blem}{\begin{lemma}}
\newcommand{\elem}{\end{lemma}}
\newcommand{\bn}{\begin{note}}
\newcommand{\en}{\end{note}}
\newcommand{\benum}{\begin{enumerate}}
\newcommand{\eenum}{\end{enumerate}}
\newcommand{\bed}{\begin{defn}}
\newcommand{\eed}{\end{defn}}
\newcommand{\br}{\begin{remark}}
\newcommand{\er}{\end{remark}}

\newtheorem{theorem}[equation]{Theorem}      
\newtheorem{lemma}[equation]{Lemma}          %
\newtheorem{corollary}[equation]{Corollary}  
\newtheorem{proposition}[equation]{Proposition}

\theoremstyle{definition}
\newtheorem{conjecture}[equation]{Conjecture}

\theoremstyle{definition}
\newtheorem{defn}[equation]{Definition}
\theoremstyle{remark}

\theoremstyle{definition}
\newtheorem{remark}[equation]{Remark}

\numberwithin{equation}{section}

\let\isom=\simeq

\let\tensor=\otimes

\newcommand{\C}{{\mathbb C}}

\newcommand{\Q}{{\mathbb Q}}

\newcommand{\Z}{{\mathbb Z}}

\renewcommand{\int}{\operatorname{int}}
\renewcommand{\O}{{\mathcal O}}

\renewcommand{\wp}{{\mathfrak p}}
\newcommand{\wq}{{\mathfrak q}}

\renewcommand{\bpro}{\begin{proposition}}
\renewcommand{\epro}{\end{proposition}}

\newcommand{\lamt}{\tilde{\lambda}}
\newcommand{\lam}{\lambda}
\newcommand{\mhalf}{M_{\half}}
\let\lab=\label
\newcommand{\half}{\frac{1}{2}}

\newcommand{\isomap}{\xymatrix{{}\ar[r]^{\isom} &}}

\newcommand{\ql}{\Q_{\ell}}

\newcommand{\gk}{G_{K}}
\newcommand{\frp}{{\rm Frob}_\wp}
\newcommand{\Tr}{{\rm Tr}}
\newcommand{\rhol}{\rho_\ell}
\newcommand{\rholp}{\rho_\ell'}

\newcommand{\vlp}{V_\ell'}

\newcommand{\vp}{V_p}
\newcommand{\vl}{V_\ell}
\newcommand{\vpp}{V_p'}
\renewcommand{\cong}{\equiv}
\newcommand{\wl}{{\mathfrak{l}}}
\begin{document}

\title[]{Ordinary reductions of abelian varieties}%
\author{Kirti Joshi}%
\address{Math. department, University of Arizona, 617 N Santa Rita, Tucson
85721-0089, USA.} \email{kirti@math.arizona.edu}

\thanks{}%


\begin{abstract}
I show that a conjecture of Joshi-Rajan on primes of Hodge-Witt reduction and in particular a conjecture of Jean-Pierre Serre on primes of good, ordinary reduction  for an abelian variety over a number field follows from a certain conjecture on Galois representations which may perhaps be easier to prove (and I prove this conjecture for abelian compatible systems of a suitable type).
This reduction (to a conjecture about certain systems of Galois representations) is based on a new slope estimate for non Hodge-Witt abelian varieties. In particular for any abelian variety over a number field with at least one prime of good ordinary or split toric reduction,  I show that the conjecture of Joshi-Rajan and the conjecture of  Serre on ordinary reductions can be reduced to proving that a certain rational trace of Frobenius is in fact an integer. The assertion that this trace is an integer is proved for abelian systems of Galois representations (of suitable type).
\end{abstract}
\maketitle
\epigraph{It don't mean a thing,\\ $\qquad$if it ain't got that swing...
}{Louis Armstrong and Duke Ellington\\
(and Irving Mills)}


\section{Introduction}
Let $k$ be a perfect field. Let $X/k$ be a smooth, projective variety over $k$. Let $W=W(k)$ be the ring of Witt vectors of $k$. Let us say, following \citep{illusie79b} that \textit{$X$ is Hodge-Witt} if the de Rham-Witt cohomology groups $H^i(X,W\Omega^j_X)$ are finitely generated as $W$-modules for all $i,j\geq0$. (see \citep{illusie79b}). Let us say that \textit{$X$ is ordinary} if and only if $H^i(X,B\Omega^j_X)=0$ for all $i,j\geq 0$, where $B\Omega^j_X=d(\Omega^{j-1}_X)$ is the sheaf of locally exact $j$-forms on $X$ (see \citep{illusie79b}). It is a theorem of \citep*{illusie83b} that if $X$ is ordinary then $X$ is Hodge-Witt. For examples of ordinary varieties see \citep{illusie90a}. 

Suppose now that $X$ is an abelian variety. Then $X$ is ordinary if and only if $X$ has $p$-rank equal to $\dim(X)$; and it is a less well-known theorem of Torsten Ekedahl (see \citep{illusie83a}) that $X$ is Hodge-Witt if and only if $p$-rank of $X$ is at least $\dim(X)-1$. In particular one sees from this that Hodge-Witt but non-ordinary abelian varieties exist.

Readers unfamiliar with \citep{illusie79b}, \citep{illusie83b}, \citep{ekedahl85,ekedahl-diagonal} should consult  Section~\ref{sec:basics} (especially Remark~\ref{r:working-def}) for a ``working definition'' of Hodge-Witt varieties adapted to abelian varieties which is more than adequate for reading this paper.

Now suppose $K/\Q$ is a  number field and  that $X/K$ is a smooth, projective variety over $K$. Then Jean-Pierre Serre has conjectured that $X$ has ordinary reduction modulo infinitely many non-archimedean primes of $K$. In \citep{joshi00a} it was conjectured that there exist infinitely many primes at which $X$ has Hodge-Witt reduction. Clearly if $X$ has good ordinary reduction at a finite prime $\wp$ of $K$ then $X$ has good, Hodge-Witt reduction at $\wp$. Thus the set of  primes of good ordinary reduction for $X$ are contained in the set of primes of good, Hodge-Witt reduction for $X$. In \citep[Theorem~4.1.3]{joshi14} I give, amongst other results, several examples which show that the two sets of primes  can have different densities.

In \citep[Theorem~4.1.3]{joshi14} I also showed that Serre's ordinarity conjecture for $X$ is equivalent to the conjecture on Hodge-Witt reductions due to Joshi-Rajan (also see \citep{joshi00a}) for $X\times_K X$ (this uses an important result of \citep{ekedahl85}) and I also proved both the conjectures for abelian varieties with complex multiplication (CM case).

In Theorem~\ref{th:hw-reduction} of this note I prove, by using methods quite different from those pursued in \citep[Theorem~4.1.3]{joshi14} for the CM case, that a certain conjecture on Galois representations (see Conjecture~\ref{galoisconj}) implies that there exist infinitely many primes of good, Hodge-Witt reductions for any  abelian variety. By \citep[Theorem~4.1.3]{joshi14} (stated here as Theorem~\ref{th:equivalence}) this is enough to prove, assuming Conjecture~\ref{galoisconj}, the result conjectured by Serre that any abelian variety over a number field has infinitely many primes of good ordinary reduction (see Theorem~\ref{th:serre}).

Let me remind the reader that for abelian surfaces the result was proved unconditionally by \citep{ogus82}. In \citep{joshi00a}, Rajan and I showed, again unconditionally, that there exist infinitely many primes of Hodge-Witt reduction for abelian threefolds. Serre's conjecture on ordinary reductions for abelian varieties has also been established in a few cases under restrictive assumptions on dimensions and endomorphism algebras or  Mumford-Tate groups  in \citep{pink98,tankeev99,noot00} and the references therein for additional results. The methods of this paper have little overlap with these works.

The proof given here proceeds in two steps. I use an observation of \citep[Theorem~4.1.3]{joshi14} (recalled here as Theorem~\ref{th:equivalence}) which allows one to reduce to proving that $A\times A$ has infinitely many primes of Hodge-Witt reduction (see Theorem~\ref{th:equivalence}). Then I prove, assuming Conjecture~\ref{galoisconj}, the existence of infinitely many primes of Hodge-Witt reduction for any  abelian variety (see Theorem~\ref{th:hw-reduction}). Finally one gets from this, assuming Conjecture~\ref{galoisconj}, the ordinarity result conjectured by Serre (see Theorem~\ref{th:serre}). The key new tool in the proofs of these two results, apart from reduction to the Hodge-Witt case provided by \citep[Theorem~4.1.3]{joshi14}, is a (sharp) slope estimate (see Theorem~\ref{hw-estimate} and Remark~\ref{rem:sharp}) for non-Hodge-Witt abelian varieties. In Theorem~\ref{th:cm-case-galoisconj}, I prove Conjecture~\ref{galoisconj} for abelian systems of Galois representations.

In  Theorem~\ref{th:toric-case-hw} I show that any abelian variety with at least one prime of good ordinary or split toric reduction over a number field has infinitely many primes of Hodge-Witt reduction and hence any such abelian variety also has infinitely many primes of ordinary reduction  provided that Conjecture~\ref{le:toric}, which reduces the proof of Conjecture~\ref{galoisconj} to proving that a certain (rational) trace of Frobenius is in fact an integer, is true. In Theorem~\ref{pro:abelian-case-effectivity} I prove Conjecture~\ref{le:toric} for abelian systems of Galois representations.


I thank Bryden Cais for many conversations and especially for patiently  listening to several of my unsuccessful  attempts to prove Conjecture~\ref{galoisconj} in all generality. I also thank Adrian Vasiu for a careful reading of an earlier version of this manuscript and a number of suggestions which led to several improvements and readability of this manuscript. Thanks are also due to Brian Conrad for some conversations on abelian Galois representations.

\section{non Hodge-Witt abelian varieties}\label{sec:basics}
Let $X$ be a smooth, projective variety over a perfect field $k$ of characteristic $p>0$. Let $W=W(k)$ be the ring of Witt vectors of $k$ and let its quotient field be $K_0$.
One says, following \citep{illusie83b}, that $X$ is Hodge-Witt if and only if the de Rham-Witt cohomology groups $H^i(X,W\Omega_X^j)$ are of finite type $W$-modules for all $i,j\geq 0$. Suppose now that $X$ is an abelian variety. It is a theorem of Torsten Ekedahl (see \citep{illusie83a}) that $X$ is Hodge-Witt if and only if $p$-rank of $X$ is $\geq \dim(X)-1$. 

\br\label{r:working-def}
Readers unfamiliar with Hodge-Witt varieties may take this $p$-rank condition as an \textit{ad hoc} definition of Hodge-Witt abelian varieties. In other words readers may adopt as a ``working definition'' the statement that $X$ is Hodge-Witt if and only if $p$-rank of $X$ is $\geq \dim(X)-1$.
Since $p$-rank of an abelian variety is additive on taking products of abelian varieties, it follows  that if $X,Y$ are abelian varieties over $k$ then $X\times Y$ is Hodge-Witt if and only if one of $X,Y$ is Hodge-Witt and the other is ordinary (this is a very special case of a beautiful general theorem of \citep{ekedahl85,ekedahl-diagonal}). This observation for abelian varieties is adequate for reading this paper. In \citep{conrad-book} abelian varieties $X$ which have $p$-rank $\geq \dim(X)-1$ are called \textit{almost ordinary} and $p$-divisible groups arising from them are said to be of \textit{extended Lubin-Tate type}. It should be noted however that neither of these labels reveal the most important property of these abelian varieties: the finiteness of the de Rham-Witt cohomology equivalently the degeneration of the slope spectral sequence at $E_1$ (see \citep{illusie83b}).
\er

In this section I prove the following slope estimate. As is conventional, a $p$-valuation $v$ used for computing slopes is normalized so that $v(p)=1$. 
\bthm\label{hw-estimate}
Let $X$ be an abelian variety over a perfect field $k$ of characteristic $p>0$. Assume that $\dim(X)=g$  and that $X$ is not Hodge-Witt. Then every slope $\lamt$ of $H^g_{\rm cris}(X/W)$ satisfies $$\lamt\geq1.$$
\ethm

\bp
Let $M=H^1_{\rm cris}(X/W)\tensor_W K_0$, for a slope $\lam$ of $M$, let $M_\lam$ be the slope $\lam$ part of $M$ and in particular let $M_0,M_1,\mhalf$ be respectively the slope zero, slope one and the slope half submodules of $M$. Let $m_\lam,m_0,m_1,m$ be their dimensions respectively. Further as $X$ is not Hodge-Witt so
\be
	m_0\leq g-2,
\ee
and one has $\dim_{K_0}(M)=2g$.
So one gets
\be
2g=\dim(M_0)+\dim(M_1)+\dim(\mhalf)+\sum_{\lam\neq 0,\frac{1}{2},1}\dim(M_\lam).
\ee
By definition of $p$-rank, the $p$-rank of $X$ is $m_0$. By duality for abelian varieties  $m_\lam=m_{1-\lam}$. Therefore one can write this as
\be\label{eq:mult-equation}
2g=2m_0+m+2\sum_{0<\lam<\half}m_\lam.
\ee
Hence one has, for every slope $0<\lam<\half$ of $M$,
\be\label{eq:nhw-bound}
2m_\lam\leq 2g-2m_0-m.
\ee
Observe that $g-m_0-\frac{m}{2}$ cannot be negative by \eqref{eq:mult-equation}.
Now the proof is split into two cases according to whether $g-m_0-\frac{m}{2}\neq 0$ or $g-m_0-\frac{m}{2}=0$.

First suppose $g-m_0-\frac{m}{2}\neq 0$. So $m_\lam \leq g-m_0-\frac{m}{2}$. Now for $0<\lam<1$ let $\lam=\frac{a}{b}$ with $a\geq 1,b>1$ and $(a,b)=1$. Then $b|m_\lam$ so $b\leq m_\lam$ and hence $\frac{1}{b}\geq \frac{1}{m_\lam}$. Hence in particular for any $0<\lam<\half$ one has
\be\label{eq:tautological-bound}
\lam\geq \frac{1}{m_\lam}.
\ee
As $m_\lam\leq  g-m_0-\frac{m}{2}$, from \eqref{eq:nhw-bound} and from \eqref{eq:tautological-bound}, for any $0<\lambda<\half$, one has the fundamental estimate:
\be\lab{eq:fun-estimate}
\lam\geq \frac{1}{m_\lambda}\geq \frac{1}{g-m_0-\frac{m}{2}}.
\ee
If $0<\lambda<\half$ is a slope, then so is $1-\lambda$ and  $1>1-\lambda>\half$ and one has
\be
1-\lambda>\lambda\geq \frac{1}{g-m_0-\frac{m}{2}},
\ee
and hence one has for all slopes $\lambda\neq 0,1,\half$ the estimate:
\be\lab{eq:fun-estimate2}
\lam\geq \frac{1}{g-m_0-\frac{m}{2}}.
\ee
Now recall that $H^g_{\rm cris}(X/W)\tensor K_0=\wedge^g H^1_{\rm cris}(X/W)=\wedge^g M$ and so any slope $\lamt$ of Frobenius on $H^g_{\rm cris}$ may be computed from slopes of $M$. Any slope $\lamt$ of $H^g_{\rm cris}$ is of the form
\be
\lamt=\lam_1+\cdots+\lam_g
\ee
for some slopes $\lam_1,\lam_2,\ldots,\lam_g$ of $M$. If any of the $\lambda_j$ occurring in this expression is equal to one, then $\lamt\geq 1$ holds trivially and there is nothing to prove. So assume $\lambda_j<1$ for all $j=1,\ldots,g$. Let $i_0$ be the number of times $\lam_j=0$ and let $i$ be the number of times $\half$ occurs in the above representation. Then
\beas\label{eq:lamt-base-estimate}
	\lamt& \geq & \frac{i}{2}+(g-i_0-i)\left(\frac{1}{g-m_0-\frac{m}{2}}\right)\\
	     &  \geq & \left(\frac{g-i_0}{g-m_0-\frac{m}{2}}\right)+i\left(\half-\frac{1}{g-m_0-\frac{m}{2}}\right).
\eeas
Note that $\half-\frac{1}{g-m_0-\frac{m}{2}}>0$. This follows from \eqref{eq:mult-equation} as
\be
	g-m_0-\frac{m}{2}=\sum_{0<\lambda<\half}m_\lambda
\ee and for slopes $0<\lambda<\half$ one has $m_\lambda\geq 3$ which gives
\be
g-m_0-\frac{m}{2}\geq 3, 
\ee
and hence
\be
\frac{1}{g-m_0-\frac{m}{2}}\leq\frac{1}{3}.
\ee
So $\half>\frac{1}{3}\geq\frac{1}{g-m_0-\frac{m}{2}}$. Thus the second term on the right hand side of \eqref{eq:lamt-base-estimate} is positive.
Now as
\be
	i_0\leq m_0
\ee
so
\be
	g-i_0\geq g-m_0\geq g-m_0-\frac{m}{2}.
\ee
Thus in the case $g-m_0-\frac{m}{2}>0$ one has
\be
	\lamt\geq \frac{g-i_0}{g-m_0-\frac{m}{2}}\geq \frac{g-m_0}{g-m_0-\frac{m}{2}}\geq 1.
\ee
Let me now address the case $g-m_0-\frac{m}{2}=0$. One argues directly using \eqref{eq:mult-equation} from which one sees that $m_\lambda=0$ for all $\lambda\neq 0,\frac{1}{2},1$. From
\be
\lamt=\lam_1+\cdots+\lam_g
\ee
assuming once again that none of the slopes in the above expression for $\lamt$ are equal to one, one has the estimate
\be
	\lamt=\frac{1}{2}(g-i_0)\geq\frac{1}{2}(g-m_0)=\frac{m}{4}.
\ee
So to prove our claim that $\lamt\geq 1$ it suffices to show that if $g-m_0-\frac{m}{2}=0$ then $m\geq 4$. If $m<4$ then  as $2|m$ this says $m=2$ or $m=0$. If $m=2$ then
\be
	g=m_0+\frac{m}{2}=m_0+\frac{2}{2}=m_0+1
\ee	
and hence $m_0=g-1$ and hence $X$ is Hodge-Witt. This  contradicts our assumption that $X$ is non Hodge-Witt. Similarly if $m=0$ then $g=m_0$ hence $X$ is ordinary but this  contradicts our assumption that $X$ is non Hodge-Witt. Thus if $g-m_0-\frac{m}{2}=0$ then $m\geq 4$, $0\leq m_0\leq g-2$ and
$\lamt=\frac{m}{4}\geq 1$.
Thus one has in all cases
\be
	\lamt\geq 1.
\ee
This completes the proof of the theorem.
\ep
\br\label{rem:sharp}
Suppose $X$ is ordinary. Then $H^g_{cris}(X/W)$ has a slope zero part of rank one. So $\lamt\geq 1$ cannot hold for all $\lamt$. Similarly if $X$ is Hodge-Witt, but not ordinary, then $X$ has $p$-rank $g-1$ hence there is a slope $\lamt=\frac{1}{2}$ for $H^g_{cris}(X/W)$. So the assumption that $X$ is non Hodge-Witt cannot be relaxed to $X$ is non-ordinary in Theorem~\ref{hw-estimate}. Moreover for every $g\geq 2$ there exists, by \citep{tate69}, an abelian variety over an algebraically closed field $k$ of characteristic $p>0$ with slopes  $\{\frac{1}{g},1-\frac{1}{g}\}$ each with multiplicity $g$ and this abelian variety is manifestly  non Hodge-Witt with exactly one slope $\lamt$ with $\lamt=1$. So the estimate is the best possible in all dimensions $g\geq 2$.
\er

\section{A conjecture about certain Galois representations}
I make the following definition. Let $K$ be a finite extension of $\Q$, and fix an algebraic closure of $K$ and let $G_K$ be the Galois group of $K$ (with respect to this algebraic closure). For any finite set of primes $S$ of $K$, and a rational prime $p$, let $S_\ell=S\cup\{\wl:\wl|\ell\}$. For a $p$-adic representation of $\gk$ and prime $\wp|p$ of $K$, I will write $V_\wp$ for the restriction of the $\gk$-representation $V_p$ to the decomposition group $D_\wp$ at $\wp$. Let $D_{cris}(V_\wp)$ be the (covariant) functor constructed in \citep{fontaine94b}. For any prime $\ell$, let $\Q_\ell(-1)$ be  the $\ell$-adic Tate twist. This is a $G_K$-representation of Weil weight two and the unique slope of Frobenius on $D_{cris}(\Q_p(-1))$ is one. With this convention $\Q_\ell(1)$ has Weil-weight minus two. 

\newcommand{\glopro}{\ref{unram}--\ref{indep0}}
\newcommand{\locpro}{\ref{pst}--\ref{indep2}}
\newcommand{\twipro}{\ref{twist}}
\newcommand{\spro}{\ref{slopes}}
\newcommand{\effpro}{\ref{integral}}
\newcommand{\hpro}{\ref{htw}}
\newcommand{\allpro}{\glopro, \locpro, \spro, \effpro, \hpro\ and \twipro}
\newcommand{\allppro}{\glopro, \locpro\ and \effpro}

\begin{defn}\label{def:twist-coupled} 
Let $g\geq 2$ be an integer. Let $S$ be a finite set of primes of $K$. A pair  of continuous, $\ell$-adic Galois representations (one for each rational prime $\ell$), $\{V_\ell\}_\ell, \{V_\ell'\}_\ell$ of $G_K$ are said to be a \textit{twist-coupled system of Galois representations of type $(g,g-2)$} if they satisfy the following conditions:
\newcommand{\Sl}{S_\ell}
\benum[label={(${\bf G}$.\arabic*)}]
\item\label{unram} Each $V_\ell$ is unramified outside $S\cup\{\ell\}$.
\item\label{weight}  Each $V_\ell$  is pure of Weil weight $g$.
\item\label{indep0} For any prime $\wp\not\in\Sl$ the characteristic polynomial of Frobenius of $V_\ell$ at $\wp$ is a polynomial $f_{\wp,\ell}(t) \in\Q[t]$, and  if  $\ell,\ell'$ are primes and   $\wp\not\in S_\ell\cup S_{\ell'}$ then $f_{\wp,\ell}(t)=f_{\wp,\ell'}(t)$.
\eenum
\benum[label={(${\bf L}$.\arabic*)}]
\item\label{pst} For all primes $\wp$, $V_\wp$ is potentially semistable.
\item\label{crys} For all primes $\wp\not\in S$, $V_\wp$  is crystalline.
\item\label{indep2} For all primes $\wp\not\in S$ the characteristic polynomial of the linearized Frobenius of $D_{cris}(V_\wp)$ is $f_\wp(t)\in\Q[t]$ and if $\wp\notin\Sl$ one has $f_{\wp,\ell}(t)=f_\wp(t)$.
\eenum

\benum[label={(${\bf E}$)}]
\item\label{integral} Except for a set of finitely many primes $\wp$ including those in $S$, the characteristic polynomial $f_\wp(t)\in\Z[t]$.
\eenum

\benum[label={(${\bf S}$)}]
\item\label{slopes}  For all primes outside the exceptional set of primes of \effpro\ the slopes of Frobenius on $D_{cris}(V_\wp)$  are in the interval $[1,g-1]$.
\eenum

\benum[label={(${\bf H}$)}]
\item\label{htw} For all primes $\wp$, $V_\wp$  has Hodge-Tate weights in $[0,g]$, and one has ${\rm gr}^0(D_{HT}(V_\wp))\neq0$.
\eenum 

\benum[label={(${\bf T}$)}]
\item\label{twist} For all $\ell$ there exists a continuous isomorphism $V_\ell\isomap V_\ell'(-1)$ of $G_K$-modules (here $\Q_\ell(-1)$ is the Tate twist).
\eenum
\end{defn}

The following simple lemma, while not essential in my proof, explicates the condition \spro.
\blem\label{le:divisibility}
Let $f(t)\in\Z[t]$ be a non-zero polynomial whose roots are $p$-Weil numbers of weight $m\geq 1$. Suppose that $L/\Q$ is a finite extension containing all the roots of $f(t)$ and suppose for any $p$-adic valuation $v_p$ of $L$, normalized so that $v_p(p)=1$, extending the standard $p$-adic valuation of $\Q_p$, one has $v_p(\alpha)\geq 1$ for any root $\alpha$  of $f(t)$. Then there is an algebraic integer $\beta$ such that $\alpha=p\beta$.
\elem
\bp 
Let $(p)=\wp_1^{a_1}\cdots\wp_r^{a_r}$ be the prime factorization of $p$ in $L$. Then the estimate $v_{\wp_i}(\alpha)\geq 1$, with the normalization $v_{\wp_i}(p)=1$, says that $\alpha\in \wp_i^{a_i}$. Hence $\alpha\in\wp_i^{a_i}$ for all $i=1,\ldots,r$. Thus $\alpha\in\cap\wp_i^{a_i}=\wp_1^{a_1}\cdots\wp_r^{a_r}=(p)$. Therefore $\alpha=p\beta$ for some algebraic integer $\beta$ as asserted.
\ep

I will write  the quantities corresponding to $V_p'$, such as characteristic polynomials of Frobenius, traces of Frobenius etc. as primed quantities: $f'_\wp(t), f_{\wp,\ell}'(t)$ etc. Hopefully there will be no confusion with derivatives (which will not be used in this paper). 

The contents of the following remarks will clarify this list of properties of twist-coupled systems and will be used in the rest of the paper. The remarks are immediate from definitions or elementary considerations.

\br\label{re:comments}
Note that the properties are indexed by global conditions  \glopro\ and by local conditions \locpro\  which members of $\{V_\ell\}_\ell$ satisfy. The assumptions \glopro\ and \ref{pst}--\ref{indep2}\ are satisfied by all geometric Galois representations.  That \glopro, \locpro\ hold in geometric contexts is the work of many mathematicians: properties \glopro\ hold by \citep{deligne74a}, properties \locpro\ first arose in the work of \citep{fontaine82a} and have now been established in geometric contexts (see \citep{colmez02} for a detailed  bibliography). That the property \ref{indep2} holds in geometric context is due to \citep{katz74}.
\benum
\item Through \twipro\ the representation  $\{V_\ell'\}_\ell$ also satisfies \glopro\  and \locpro. 
\item Properties \effpro, \spro, \hpro\ are not invariant under Tate twists and hence are not inherited by $\{V_\ell'\}$ through \twipro.
\item Property \effpro\ of $V_p$ says that in some sense $V_p$ is ``effective.'' 
\item Note that properties \ref{slopes}, \ref{integral}\ and \ref{htw}\  only involve $\{V_\ell\}_\ell$.
\eenum
\er

\bpro\label{integralp}
Suppose $\{V_\ell,V_\ell' \}_\ell$ is a twist-coupled family of Galois representations. Then
\benum
\item Properties \glopro, \locpro\   are hold for $\{V_\ell'\}$.
\item The property \effpro\  is inherited by $\{V_\ell'\}$.
\item For  all primes $\wp\not\in S$, $V_\wp'$ is crystalline.
\item For all but a finite number of primes, the slopes of Frobenius on $D_{cris}(V_\wp')$ are in $[0,g-2]$.
\item The family  $\{V_\ell^{ss},V_\ell'^{ss} \}_\ell$ of semisimplifications of  $\{V_\ell,V_\ell' \}_\ell$ is also a twist coupled system.
\eenum
\epro

\bp 
(1), (3) and (4) are immediate from the fact that \allpro\ hold for $\{V_\ell\}$ and as $V_\ell'=V_\ell(1)$ by \twipro. Claim (2) is now immediate from \ref{twist}, \ref{integral} \ and \ref{slopes}\ for $\{V_\ell\}$. Clearly (5) is entirely formal.
\ep

\br 
In particular one sees from Proposition~\ref{integralp}(2) that if \allpro\ hold then \ref{integral} also holds for $\{V_\ell'\}_\ell$. As $V_\ell=V_\ell'(-1)$ by \ref{twist},  and as  \ref{integral} holds for $V_\ell'$, one has \emph{divisibility of traces of Frobenius elements acting on $V_\ell$}. Thus the conditions \allpro\ provide a natural way of encoding the divisibility by $p$ of trace of Frobenius element at $\wp$ acting on $V_\ell$. 
\er


I propose the following conjecture.

\bcon\label{galoisconj}
Let $X/K$ be a smooth, projective variety over a number field $K$ of dimension $g\geq 2$ with $$H^g(X,\O_X)\neq0.$$ Then the family of $\ell$-adic Galois representations:
\be
\{V_\ell=H^g_{et}(X,\Q_\ell),V_\ell'=V_\ell(1)\}_\ell
\ee
is  not  a  twist coupled system of Galois representations of type $(g,g-2)$.
\econ

For $g=2$ and any smooth projective variety $X$ (with $H^2(X,\O_X)\neq0$), Conjecture~\ref{galoisconj} is immediate by using methods of \citep{ogus82,joshi00a,bogomolov09}. To see this it is sufficient to note that the methods of \textit{loc. cit.} show that if the traces of Frobenius at all but finite number of  primes $p$  (for $\ell$-adic \'etale cohomology, with an appropriate $\ell$ for $X$) are divisible by $p$, then the  semisimplification of $H^2_{et}(X,\Q_\ell)$ is isomorphic to $\oplus\Q_\ell(-1)$ (as a $\gk$-module) (and by compatibility (\glopro) this must hold for all $\ell$) which obviously contradicts \hpro\ (by the Hodge-Tate decomposition Theorem furnished by $p$-adic Hodge Theory).  Thus $\{V_\ell=H^2_{et}(X,\Q_\ell),V_\ell'=V_\ell(1)\}_\ell$ is not a twist-coupled system.  Note that the hypothesis in \textit{loc. cit.} for $g=2$, that $X$ is  a K3 surface is used only to prove that there is exactly one unit eigenvalue of Frobenius. Here is some additional evidence for Conjecture~\ref{galoisconj}.

\bthm\label{th:cm-case-galoisconj}
Let $X$ be a smooth, projective variety of dimension $g$ over a number field and suppose $H^g(X,\O_X)\neq 0$. Assume that  $\{(\rhol, H^g_{et}(X,\ql))\}_\ell$ is a family of abelian Galois  representations (i.e. factors as a representation $\rhol:\gk^{ab}\to {\rm GL}(H^g_{et}(X,\ql))$. Then Conjecture~\ref{galoisconj} holds for $X$.
\ethm
\bp 
As $\{H_{et}^g(X,\ql)\}_\ell$ is a family of abelian representations, so is its Tate-twist $\{V_\ell'=H_{et}^g(X,\ql)(1)\}_\ell$. Replacing these representations by their semi-simplifications one can assume that one has a family of semi-simple abelian representations. Assume the Conjecture~\ref{galoisconj} is false for $X$. This means that 
\be
	\{V_\ell=H_{et}^g(X,\ql),H_{et}^g(X,\ql)(1)\}_\ell
\ee 
forms a twist-coupled system of $\gk$-representations (so properties \allpro\ hold). In particular by Proposition~\ref{integralp}(2), $\{V_\ell'\}$ satisfies \effpro. By \glopro, \locpro\ and \citep[Theorem of Tate, page {III-7}]{serre-abelian} $V_p'$ is locally algebraic and hence by \citep[Proposition, page {III-9}]{serre-abelian} $V_p'$ arises from a $\Q$-rational representation of an algebraic torus.  It is possible to choose a prime $p$ such that for any prime $\wp$ of $K$ lying over $p$ such that $V_p'$ is crystalline at $\wp$ (and hence Hodge-Tate at  $\wp|p$) and the algebraic torus giving rise to this representation is split at $p$.  Pick such a prime $p$.  Then $V_\wp'$ is a direct sum of crystalline characters of $D_\wp$.  As $\{V_\ell'\}$ satisfies \ref{integral},  the characteristic polynomials of Frobenius elements at all but finite number of primes of $K$ (acting on $V_p'$) have integer coefficients.  Then the integrality of the characteristic polynomials of Frobenius elements in $\gk$ (over all but finitely many primes of $K$) shows  by \citep[Chapter II, Corollary 2, page {II-36}]{serre-abelian} that the Hodge-Tate weights of $V_\wp'$ are non-negative (note that \textit{loc. cit.} the Hodge-Tate weights are determined by the algebraic characters of the aforementioned torus which appear in rational representation of this torus provided by \citep[Proposition, page {III-9}]{serre-abelian}).  So $V_p'$ has non-negative Hodge-Tate weights and hence $V_p=V_p'(-1)$ has strictly positive Hodge-Tate weights. This contradicts \ref{htw} as a Theorem of Tate (see \citep{tate67}) says that if $V_p$ is a $p$-adic representation with strictly positive Hodge-Tate weights (at primes over $p$) then $(V_p\tensor\C_p)^{D_\wp}=0$ for any prime $\wp$ lying over $p$. On the other hand one has $H^g(X,\O_X)\neq 0$ and hence $(V_p\tensor\C_p)^{D_\wp}\not=0$. Hence one has arrived at a contradiction.
\ep 

In some situations one can replace Conjecture~\ref{galoisconj} by the somewhat simpler Conjecture~\ref{le:toric} given below.

\bcon\label{le:toric} Let $K$ be a number field. Let $S$ be a finite set of primes of $K$.  Suppose $\{V_\ell' \}_\ell$ is a family  of continuous Galois representations satisfying the properties \glopro, \locpro, \effpro\  and suppose that for a prime $\wp$ (lying over a rational prime $p$) the following hypothesis hold (here $D_\wp\subset\gk$ is the decomposition group at $\wp$ and $I_\wp$ is the inertia subgroup at $\wp$, and $V_\wp'$ is the restriction of $V_p'$ to $D_\wp$):
\benum[label={(${\bf O}$.\arabic*)}]
\item\label{uniformity} The vector of Hodge-Tate weights of $V_\wq'$ is constant as $\wq$ varies over primes of $K$.  
\item\label{unipotence-at-toric-prime} For any prime $\ell$ not lying below $\wp$, the representation $V_\ell'$ of $D_\wp$ is unramified (resp. unipotent). 
\item\label{ordinarity-at-toric-prime} The representation $V_\wp'$ of $D_\wp$ is crystalline ordinary (resp. semi-stable ordinary, i.e., equipped with a $D_\wp$-invariant filtration whose graded pieces are isomorphic to $\Q_p(-i)$ for $i\in\Z$).
\item\label{compatibility-at-toric-prime} For any prime $\ell$ not lying below $\wp$, the characteristic polynomial of Frobenius at $\wp$ acting on $V_\ell'$ coincides with the characteristic polynomial of Frobenius on $D_{st}(V_\wp')$.
\item\label{integrality-at-toric-prime} The trace $a_\wp'$ of Frobenius $\phi_\wp'$ on $D_{st}(V_\wp')$ is rational (i.e. $a_\wp'\in\Q$).
\eenum
Then the trace $a_\wp'$  of Frobenius $\phi'_\wp$ on $D_{st}(V_\wp')$ is  an integer.
\econ

Here is some partial evidence for Conjecture~\ref{le:toric}.

\bthm\label{pro:abelian-case-effectivity}
Suppose $\{\rholp:\gk^{ab}\to{\rm GL}(V_\ell')\}_\ell$ is an  family of continuous, abelian Galois representations which satisfies all the hypothesis of Conjecture~\ref{le:toric}. 
Then Conjecture~\ref{le:toric} is true for $\{V_\ell'\}_\ell$.
\ethm
\bp 
This is proved in a manner similar to Theorem~\ref{th:cm-case-galoisconj}. Replacing our system by its semisimplification one can assume that one has a semisimple system of representations of $\gk^{ab}$. First choose a prime $q\neq p$ and a prime $\wq$ of $K$ lying over $q$ such that $V_\wq'$ is crystalline at $\wq$.  This means, by \glopro, \locpro\ and \citep[Theorem of Tate, page {III-7}]{serre-abelian}, that $V_\wq'$ is locally algebraic and hence by \citep[Proposition, page {III-9}]{serre-abelian} $V_\wq'$ arises from a $\Q$-rational representation of an algebraic torus.   As $\{V_\ell'\}$ satisfies \ref{integral},  the characteristic polynomials of Frobenius elements at all but finite number of primes (acting on $V_q'$) are integers. The integrality of the characteristic polynomials of Frobenius elements in $\gk$ (over all but finitely many primes of $K$) shows,  by \citep[Chapter II, Corollary 2, page {II-36}]{serre-abelian}, that the Hodge-Tate weights of $V_\wq'$ are non-negative. By \ref{uniformity} the Hodge-Tate weights of $V_\wp'$ are also non-negative. Since $V_\wp'$ is semistable and ordinary (by \ref{ordinarity-at-toric-prime}), the non-negativity of Hodge-Tate weights of $V_\wp'$ says that the characteristic polynomial of $\phi_\wp'$ on $D_{st}(V_\wp')$ has $p$-adic integer coefficients. Thus by \ref{integrality-at-toric-prime} one sees that the trace of Frobenius $\phi'_\wp$ is a rational integer. This proves the assertion.
\ep

\section{Hodge-Witt and Ordinary reductions}
Before proceeding further let me recall the following  observation of \citep[Theorem~4.1.3]{joshi14}. The only case Theorem~\ref{th:equivalence} of interest for this paper is the case when $X$ is an abelian variety and as remarked in Remark \ref{r:working-def} if one uses the ``working definition'' that an abelian variety $X$ is Hodge-Witt if and only if $X$ has $p$-rank $\geq\dim(X)-1$,  then Theorem~\ref{th:equivalence} given below is quite elementary to prove. 
\bthm\label{th:equivalence}
Let $X$ be a smooth, projective variety over a number field $K$. Then the following are equivalent:
\begin{enumerate}
\item There exist infinitely many primes of ordinary reduction for $X$.
\item There exist infinitely many primes of ordinary reduction for $X\times_K X$
\item There exist infinitely many primes of Hodge-Witt reduction for $X\times_K X$.
\end{enumerate}
\ethm

For $\dim(X)=3$ the existence of infinitely many primes of Hodge-Witt reductions, (without assuming Conjecture~\ref{galoisconj}), is due to \citep{joshi00a} (also see \citep{joshi14}).

\bthm\label{th:hw-reduction}
Let $K$ be a number field. Let $X/K$ be an  abelian variety over $K$ of dimension $g\geq 1$. Assume Conjecture~\ref{galoisconj} is true for $X$ if $g\geq 2$. Then there exist infinitely many primes of Hodge-Witt reduction for $X$.
\ethm

\bp
Since any abelian variety of dimension $g=1$ has Hodge-Witt for an infinite set of primes, one can assume that $g\geq 2$. Suppose the assertion is not true. Then I show that 
\be 
	\{V_\ell=H^g_{et}(X,\Q_\ell),V_\ell'=V_\ell(1)\}_\ell
\ee
is a system of twist-coupled  Galois representations of type $(g,g-2)$. Note that \twipro\ holds as one has the tautological isomorphism $V_\ell=V_\ell(1)(-1)=V_\ell'(-1)$. Clearly $\{V_\ell,V_\ell'\}_\ell$ satisfies all the properties \allpro\ except possibly \ref{slopes}. By our assumption $X$ does not have Hodge-Witt reduction at all but a finite number of primes $\wp$ of $K$. Suppose $\wp$ is a prime of good, non Hodge-Witt reduction.  Let $X_\wp$ be the reduction modulo $\wp$ of $X$ (in a regular, smooth, proper model of $X$ over a suitable localization of the ring of integers of $K$). Let $\kappa(\wp)$ be the residue field of $\wp$. Then by Theorem~\ref{hw-estimate} one sees that if $\alpha$ is any eigenvalue of Frobenius on $H^g_{et}(X,\Q_\ell)$ (equivalently on $H^g_{cris}(X_\wp/W(\kappa(\wp)))$, by \citep{katz74}) then $\alpha$ satisfies $v_p(\alpha)\geq 1$. As $X_\wp$ has dimension $g$, one has by Poincar\'e duality for crystalline cohomology of $X_\wp$ that $v_p(\alpha)\leq g-1$. Hence for any prime $\wp$ of non Hodge-Witt reduction \spro\ holds for $D_{cris}(V_\wp)$. Thus one has arrived at a twist coupled system of Galois representations of type $(g,g-2)$. So if  Conjecture~\ref{galoisconj} is true for $X$ then one has arrived at a contradiction.
\ep

\bcor\label{cor:density}
Let $X/K$ be an abelian variety over a number field $K$. If Conjecture~\ref{galoisconj} is true for $X$ then there exists a set of primes of positive density at which $X$ has Hodge-Witt reduction.
\ecor

\bp
If the assertion is not true then the set of primes of  Hodge-Witt reduction has density zero and examining the proof of Theorem~\ref{th:hw-reduction} one is again  led to a contradiction.
\ep

\br
Let me point out that if one assumes that $X$ has non-ordinary reduction at $\wp$ in Theorem~\ref{hw-estimate} then it follows that the traces of Frobenius of $\wp$ on $H^g_{et}(X,\Q_\ell)$ are divisible by $p$. However this is not enough to conclude that all eigenvalues of Frobenius are divisible by $p$. Thus the assumption in Theorem~\ref{hw-estimate} that $X$ is non Hodge-Witt  places stronger constraints on the crystalline cohomology of $X$ than non-ordinarity does.
\er

\bthm[Serre's Ordinarity Conjecture]\label{th:serre}
Let $A/K$ be any abelian variety over $K$. If Conjecture~\ref{galoisconj} is true  for $A\times_K A$  then there exists a set of  primes of positive density of good ordinary reduction for $A$.
\ethm
\bp
It is sufficient, by \citep[Theorem~4.1.3]{joshi14} (also see Remark~\ref{r:working-def}), to prove that $X=A\times_K A$ has infinitely many primes of Hodge-Witt reduction. This is immediate from Theorem~\ref{th:hw-reduction} if $A\times_KA$ satisfies Conjecture~\ref{galoisconj}. The density assertion follows from Theorem~\ref{th:hw-reduction} and Corollary~\ref{cor:density}.
\ep

\bcor
Let $K$ be any field of characteristic zero. If Conjecture~\ref{galoisconj} is true for an abelian variety and also for the self-product of this variety over $K$ then this abelian variety  has infinitely many primes of Hodge-Witt and ordinary reductions.
\ecor
\bp
This is immediate from Theorem~\ref{th:serre}.
\ep

\br 
In \citep{joshi14} it is conjectured that for a class of smooth, projective varieties including abelian varieties, the set of primes of non Hodge-Witt reductions in dimension $\dim(X)$ is also infinite. This, I believe, is the correct analog in higher dimensions, of Elkies' Theorem (see \citep{elkies87}) on the infinitude of primes of supersingular reductions for elliptic curves over $\Q$ (in \citep{joshi14} it is noted that Elkies's result is a very special case of this conjecture).
\er 

\section{The Toric Case}\label{se:toric-case}
In case of abelian varieties over a number field with at least one prime of ordinary or split-toric reduction one may replace Conjecture~\ref{galoisconj} by  Conjecture~\ref{le:toric} in the proof of Theorem~\ref{th:hw-reduction}. This is proved in Theorem~\ref{th:toric-case-hw} below.
 
\bthm\label{th:toric-case-hw}
Let $X/K$ be an abelian variety of dimension $g\geq 2$ over a number field $K$ and suppose that $\wp$ is a prime of $K$ at which $X$ has either ordinary or split toric reduction. If Conjecture~\ref{le:toric} is true then  $X$ has infinitely many primes of Hodge-Witt reduction.
\ethm
\bp 
Let us suppose the theorem is false.  Then the proof of Theorem~\ref{th:hw-reduction} shows that $\{V_\ell=H^g_{et}(X,\Q_\ell),V_\ell'=V_\ell(1) \}_\ell$  forms a twist-coupled system of Galois representations. So Properties \allpro\ hold.

Let $a_\wp$ (resp. $a_\wp'$) be the trace of Frobenius on $D_{st}(\vp)$ (resp. $D_{st}(\vpp)$). Suppose for the moment that $a_\wp$ is an integer. Then the property \twipro\ says that $a_\wp'$ is rational. Since $\vp$ is ordinary and $H^g(X,\O_X)\neq0$ one has
\be\label{eq:indiv}
a_\wp\not\cong0\bmod{p}.
\ee
On the other hand the isomorphism $\vl\isom\vlp(-1)$ of $\gk$-modules says that
one has 
\be\label{eq:trace-rel}
a_\wp=a_\wp'\cdot p.
\ee 

Now  it is clear that \eqref{eq:indiv} and \eqref{eq:trace-rel} cannot hold simultaneously provided one knows that $a_\wp'\in\Z$. 

First suppose that $X$ has good ordinary reduction at $\wp$. Suppose $\ell$ is not a prime lying below $\wp$. Then by the hypothesis, that $X$ has good reduction at $\wp$, the trace of Frobenius at $\wp$ acting on $V_\ell$ is an integer. By \ref{indep2} this integer is also the trace of Frobenius acting on $D_{st}(V_\wp)=D_{crys}(V_\wp)$. As $V_\wp'=V_\wp(1)$ this says that the trace of Frobenius at $\wp$ acting on $V_\ell'$ is rational. Hence all the hypothesis of Conjecture~\ref{le:toric} hold for $\{\vlp\}_\ell$. Thus one can invoke Conjecture~\ref{le:toric} and the assertion follows.

Now suppose that $\wp$ is a prime of split toric reduction lying over a rational prime $p$. Let us consider the restriction of the $G_K$-representation of $V_\ell$ (for a rational prime $\ell$ not lying below $\wp$) to the decomposition group $D_\wp$ at $\wp$.  Let $K_\wp$ be the completion of $K$ at $\wp$.

I have to prove that the traces $a_\wp$ (resp. $a_\wp'$) are rational. This is certainly well-known but I recall this for convenience.

Let  $\sigma$ be any element of the Weil group of $K_\wp$ (over its algebraic closure) lifting  Frobenius of the residue field at $\wp$ and  acting on $\vl$ (through the action of $D_\wp$ on $\vl$). By abuse of terminology I will call such an element a Frobenius at $\wp$. Then it follows from \citep{sga7-1} and \citep{coleman99} that $a_\wp=\Tr(\rhol(\frp))$ and it is immediate from \citep[Theorem~4.3(b), page 359]{sga7-1}) and the fact that as $D_\wp$-modules one has an isomorphism $V_\ell=H^g_{et}(X,\Q_\ell)=\wedge^g H^1_{et}(X,\Q_\ell)$  that $a_\wp$ is an integer and so $a_\wp'$ is rational. Hence all the hypothesis of Conjecture~\ref{le:toric} hold for $\{\vlp\}_\ell$.
Now one invokes Conjecture~\ref{le:toric} to see that $a_\wp'$ is an integer. This together with \eqref{eq:trace-rel} and \eqref{eq:indiv} completes the proof.
\ep

This has the following corollary.
\bthm~\label{th:toric-case}
Let $X/K$ be an abelian variety over a number field $K$ and suppose that $\wp$ is a prime of $K$ at which $X$ has good ordinary or split toric reduction. Assume that Conjecture~\ref{le:toric} is true. Then $X$ has infinitely many primes of ordinary reduction.
\ethm
\bp 
This follows from Theorem~\ref{th:toric-case-hw} and Theorem~\ref{th:equivalence} because if $X$ has a prime of split toric reduction then so does $X\times X$. Similarly if $X$ has one prime of good ordinary reduction then so does $X\times_KX$. Now the assertion follows from Conjecture~\ref{le:toric}.
\ep

\bibliographystyle{plainnat}
\bibliography{serre,../../master/master6}

\begin{thebibliography}{25}
\providecommand{\natexlab}[1]{#1}
\providecommand{\url}[1]{\texttt{#1}}
\expandafter\ifx\csname urlstyle\endcsname\relax
  \providecommand{\doi}[1]{doi: #1}\else
  \providecommand{\doi}{doi: \begingroup \urlstyle{rm}\Url}\fi

\bibitem[Bogomolov and Zarhin(2009)]{bogomolov09}
Fedor Bogomolov and Yuri Zarhin.
\newblock Ordinary reduction of {K3} surfaces.
\newblock \emph{Cent. Eur. J. Math.}, 7\penalty0 (2):\penalty0 206--213, 2009.

\bibitem[Chai et~al.(2014)Chai, Conrad, and Oort]{conrad-book}
Ching-Li Chai, Brian Conrad, and Frans Oort.
\newblock \emph{Complex multiplication and lifting problems}, volume 197 of
  \emph{Mathematical Surveys and Monographs}.
\newblock American {M}ath. {S}ociety, 2014.

\bibitem[Coleman and Iovita(2009)]{coleman99}
Robert Coleman and Adrian Iovita.
\newblock The frobenius and the monodromy operator for curves and abelian
  varieties.
\newblock \emph{Duke Math. Journal}, 97\penalty0 (1):\penalty0 171--215, 2009.

\bibitem[Colmez(2001-2002)]{colmez02}
Pierre Colmez.
\newblock Les conjectures de monodromie $p$-adiques.
\newblock \emph{S\'eminaire {N}. {B}ourbaki}, exp. 897, 2001-2002.

\bibitem[Deligne(1974)]{deligne74a}
P.~Deligne.
\newblock La conjecture de {Weil} {I}.
\newblock \emph{Publ. Math. I.H.E.S}, 43:\penalty0 5--77, 1974.

\bibitem[Ekedahl(1985)]{ekedahl85}
Torsten Ekedahl.
\newblock On the multiplicative properties of the de {Rham}-{Witt} complex
  {II}.
\newblock \emph{Ark. f\"ur {Mat}.}, 23, 1985.

\bibitem[Ekedahl(1986)]{ekedahl-diagonal}
Torsten Ekedahl.
\newblock \emph{Diagonal complexes and {$F$}-guage structures}.
\newblock Travaux ex {C}ours. Hermann, Paris, 1986.

\bibitem[Elkies(1987)]{elkies87}
Noam Elkies.
\newblock The existence of infinitely many supersingular primes for every
  elliptic curve over ${\Q}$.
\newblock \emph{Invent. {M}ath.}, 89\penalty0 (3):\penalty0 561--567, 1987.

\bibitem[Fontaine(1982)]{fontaine82a}
Jean-Marc Fontaine.
\newblock Sur certains types de repr\'esentations $p$-adiques du groupe de
  galois d'un corps local; construction d'un anneau de barsotti-tate.
\newblock \emph{Annals of Maths.}, 115:\penalty0 529--572, 1982.

\bibitem[Fontaine(1994)]{fontaine94b}
{J}ean-{M}arc Fontaine.
\newblock Expose {III}: R\'epresentations $p$-adiques semi-stables.
\newblock In \emph{$p$-adic Periods}. 1994.

\bibitem[Grothendieck(1972)]{sga7-1}
A.~Grothendieck.
\newblock \emph{Groupes de monodromie en g\'eom\'etrie alg\'ebrique,
  S\'eminaire de Geometrie algebrique du Bois-Marie 1967-69}.
\newblock Number 288 in Springer Lecture Notes in Math. Springer-Verlag, 1972.

\bibitem[Illusie(1979)]{illusie79b}
Luc Illusie.
\newblock Complexe de de {R}ham-{W}itt et cohomologie cristalline.
\newblock \emph{Ann. Scient. Ecole Norm. Sup.}, 12:\penalty0 501--661, 1979.

\bibitem[Illusie(1983)]{illusie83a}
Luc Illusie.
\newblock \emph{Algebraic {G}eometry {Tokyo/Kyoto}}, volume 1016 of
  \emph{Lecture {Notes} in {M}athematics}, chapter Finiteness, duality and
  K\"unneth theorems in the cohomology of the de {Rham}-{Witt} complex, pages
  20--72.
\newblock Springer-{Verlag}, 1983.

\bibitem[Illusie(1990)]{illusie90a}
Luc Illusie.
\newblock \emph{The {G}rothendieck {F}estschrift, Vol 2.}, volume~87 of
  \emph{Progr. {M}ath.}, chapter Ordinarit\'e des intersections compl\`etes
  g\'en\'erales, pages 376--405.
\newblock Bikhauser, 1990.

\bibitem[Illusie and {R}aynaud(1983)]{illusie83b}
Luc Illusie and Michel {R}aynaud.
\newblock Les suites spectrales associ\'ees au complexe de de {R}ham-{W}itt.
\newblock \emph{Inst. {H}autes {\'E}tudes {S}ci. {P}ubl. {M}ath.}, 57:\penalty0
  73--212, 1983.

\bibitem[Joshi(2014)]{joshi14}
Kirti Joshi.
\newblock On primes of ordinary and hodge-witt reductions.
\newblock Preprint: \url{}, 2014.

\bibitem[Joshi and Rajan(2000)]{joshi00a}
Kirti Joshi and C.~S. Rajan.
\newblock On frobenius split and ordinary varieties.
\newblock \emph{Preprint: \url{http://arxiv.org/abs/math/0110070}}, 2000.

\bibitem[Katz and Messing(1974)]{katz74}
Nicholas Katz and William Messing.
\newblock Some consequences of the {R}iemann hypothesis for varieties over
  finite fields.
\newblock \emph{Invent. {M}ath.}, 23:\penalty0 73--77, 1974.

\bibitem[Noot(2000)]{noot00}
Rutger Noot.
\newblock Abelian varieties with $\ell$-adic {G}alois representations of
  {M}umford's type.
\newblock \emph{Journal f\"ur die reine und angewandte Mathematik 519 (), p.
  .}, \penalty0 (number), 2000.

\bibitem[Ogus(1982)]{ogus82}
Arthur Ogus.
\newblock \emph{Hodge cycles, {M}otives, {S}himura varieties}, volume 900 of
  \emph{Springer {L}ecture {N}otes in {M}athematics}, chapter Hodge cycles and
  crystalline cohomology.
\newblock Springer-{V}erlag, 1982.

\bibitem[Pink(1998)]{pink98}
Richard Pink.
\newblock $\ell$-adic algebraic monodromy groups, cocharacters , and the
  {M}umford-{T}ate conjecture.
\newblock \emph{J.~{R}eine {A}ngew. {M}ath.}, 495:\penalty0 187--237, 1998.

\bibitem[Serre(1968)]{serre-abelian}
J.-P. Serre.
\newblock \emph{Abelian $\ell$-adic representations}.
\newblock Benjamin, New {Y}ork, 1968.

\bibitem[Tankeev(1999)]{tankeev99}
S.~G. Tankeev.
\newblock On weights of the $\ell$-adic representation and arithmetic of
  frobenius eigenvalues.
\newblock \emph{Izvestiya {M}ath.}, 63\penalty0 (1):\penalty0 181--218, 1999.

\bibitem[Tate(1967)]{tate67}
J.~Tate.
\newblock $p$-divisible groups.
\newblock In \emph{Proc. conf. on Local fields at Driebergen}, pages 158--183,
  Berlin, 1967. Springer.

\bibitem[Tate(1968-1969)]{tate69}
John Tate.
\newblock Classes d'isog\'enie des vari\'et\'es ab\'eliennes sur un corps fini.
\newblock \emph{S\'eminaire {N}. {B}ourbaki}, exp. 352:\penalty0 95--110,
  1968-1969.

\end{thebibliography}
\end{document}